\documentclass{ifacconf}

\usepackage{graphicx}      
\usepackage{natbib}        
\usepackage{amsmath,amsfonts,amssymb,amscd}
\renewcommand{\qed}{\hfill \mbox{\raggedright \rule{0.1in}{0.1in}}}
\newcommand{\R}{\mathbb{R}}
\newcommand{\Z}{\mathbb{Z}}
\newcommand{\trp}{\mathrm{T}}
\def\ker{\mathop{\rm ker}\nolimits}
\def\ran{\mathop{\rm ran}\nolimits}
\def\rank{\mathop{\rm rank}\nolimits}
\def\sp{\mathop{\rm span}\nolimits}

\begin{document}
\begin{frontmatter}

\title{
On the Realization of Hidden Markov Models and Tensor Decomposition
} 

\author[First]{Yoshito Ohta}

\address[First]{Kyoto University, 
   Kyoto, 606-8501 Japan (e-mail: yoshito\_ohta@i.kyoto-u.ac.jp).}

\begin{abstract}                
The minimum realization problem of hidden Markov models (HMM's) is a
 fundamental question of stationary discrete-time processes with a
 finite alphabet.  It was shown in the literature that tensor
 decomposition methods give the hidden Markov model with the minimum
 number of states generically.  However, the tensor decomposition
 approach does not solve the minimum HMM realization problem when the
 observation is a deterministic function of the state, which is an
 important class of HMM's not captured by a generic argument.  In this
 paper, we show that the reduction of the number of rank-one tensors
 necessary to decompose the third-order tensor constructed from the
 probabilities of the process is possible when the reachable subspace is
 not the whole space or the null space is not the zero space.  In fact,
 the rank of the tensor is not greater than the dimension of the
 effective subspace or the rank of the generalized Hankel matrix.
\end{abstract}


\end{frontmatter}

\section{Introduction}
\label{sec:introduction}

A hidden Markov model (HMM) produces a finite-valued process as the
output of a finite-state Markov process.  Because of the ability to
model various kinds of signals, HMM's have been exploited to solve many
real-world problems such as speech processing 
(\cite{Rabiner1989}) and computational biology (\cite{Krough1994}). 

The realization problem of HMM's is to derive a finite-state Markov
model with an observation map 
given the statistics of the stochastic process.  
This direction of research started in 
\cite{BlackwellKoopmans1957} 
and
\cite{Gilbert1959}, where they considered whether 
a finite-state Markov model could be 
uniquely identified from the statistics of the output process.  
The studies by 
\cite{Anderson1999} and 
\cite{Vidyasagar2011} provide 
a comprehensive overview and original results on the realization problem.  

As is discussed in \cite{Anderson1999}, the realization problem is
closely related to the so-called generalized Hankel matrix whose
entries are probabilities of strings occurring in the process arranged
in a certain way just as the conventional Hankel matrix formed from 
impulse response coefficients.  
The finite-rank property of the generalized Hankel matrix is a necessary
condition for the steady-state process to have a finite-state HMM
realization. Still, the converse does not hold in general.  
See the discussion and the example in \cite{Vidyasagar2011}.  
The difficulty of characterizing the existence of a finite-state Markov
model is the non-negativity of the transition matrix.  
A realization where the non-negativity constraint is relaxed is called 
a quasi-realization or a pseudo-realization and studied in 
\cite{ItoAmariKobayashi1992,Anderson1999,Vidyasagar2011}. 

Recently, \cite{HuangGeKakadeDahleh2016} studied the minimal realization
problems for HMM's and showed that a minimal quasi-HMM-realization and a
minimal HMM-realization could be efficiently solved 
generically for almost all HMM's. 
The main tool for the calculation of a minimal HMM-realization is
third-order tensor decomposition. 
It was shown that the rank of the tensor is equal to the minimum degree
of an HMM realization using a sufficient condition for 
the uniqueness of tensor decomposition studied in \cite{Kruskal1977} 
(see also the survay paper by \cite{KoldaBader2009}). 
The paper assumes that any two columns of the observation matrix are
linearly independent.  If this is not the case, there are two identical
columns.  However, unlike the claim in \cite{HuangGeKakadeDahleh2016}, 
two states cannot be merged to give an equivalent HMM
realization of small order because the HMM-realization problem is
nontrivial even if the observation is a deterministic function of the
state. 

In this note, we consider the uniqueness of third-order tensor
decomposition when the observation is a deterministic function of the
state and show that the tensor has low-rank decomposition if the
generalized Hankel matrix has lower rank than the number of states.  
The study by \cite{ItoAmariKobayashi1992} introduced the notion of the
reachable space and the null space of HMM's. 
Notice that the condition implies that 
either the reachable space is not the whole space or 
the null space is nontrivial. 
Hence the tensor decomposition approach does not solve the minimum
HMM realization problem when the observation is a deterministic function
of the state, which is an important class of HMM's not captured
generically. 

The rest of the paper is organized as follows. 
In Section \ref{sec:realization_problem}, the realization problem of 
HMM's is described when a stationary process is given.  
Section \ref{sec:reachable_and_null_spaces} reviews
the notion of the reachable space and the null space of HMM's  
and derives representations of these subspaces.  
Section \ref{sec:tensor_decomposition} considers rank reduction of
tensor decomposition when the reachable space is not the whole space, and
the null space is not the zero subspace.  

We use the following notations; 
$\Z$ denotes the set of integers, 
$\R$ denotes the set of real numbers, 
$\R^{m}$ denotes the set of real vectors 
of size $m$,  
and 
$\R^{m \times n}$ denotes the set of real matrices 
of size $m \times n$. 
For a matrix $A \in \R^{m \times n}$ or a vector $a \in \R^{m}$, 
$A^{\trp}$ or $a^{\trp}$ denotes the transposition of the matrix or the
vector. If $A \in \R^{m \times m}$ is invertible, $A^{-\trp}$ denotes the
inverse of $A^{\trp}$.

\section{Realization Problem}
\label{sec:realization_problem}

Suppose $\left\{ y_{t} \right\}$ is a stationary discrete-time random
process taking values in a finite set 
$\left\{ 1, \ldots, d \right\}$.  
For $t,s \in \Z$, 
let $u_{s}^{t} = \left( u_{s}, \ldots, u_{t} \right) \in 
\left\{ 1,\ldots,d \right\}^{\left| t-s \right|+1}$ is an array of
length $\left| t-s \right|+1$. Note that the array is in ascending order
if $s<t$ and in descending order if $s>t$. 
Let the random vector $y_{s}^{t} = \left( y_{s}, y_{s+1},
\ldots, y_{t} \right)$ be defined similarly.  
Let 
$\mathbb{P}(y_{s}^{t}=u_{s}^{t})$ denote the probability of the event
$y_{s}^{t}= u_{s}^{t}$.  
From the stationarity, 
$\mathbb{P}(y_{s}^{t}=u_{s}^{t})=
\mathbb{P}(y_{s+\tau}^{t+\tau}=u_{s}^{t})$ for any 
$s,t,\tau \in \Z$ and 
$u_{s}^{t} \in \left\{ 1,\ldots,d \right\}^{\left| t-s \right|+1}$. 

Suppose $\left\{ x_{t} \right\}$ is a stationary Markov chain taking
values in a state space $\left\{ 1, \ldots, k \right\}$ with the state 
transition matrix $Q = \left( q_{ij} \right) \in \R^{k \times k}$ such
that
\[
\mathbb{P} \left( x_{t+1} = i \mid x_{t} = j \right) = q_{ij}, \quad
i,j \in \left\{ 1,\ldots,k \right\}.
\]
Note that $Q$ is nonnegative and column-stochastic, i.e., the sum of every
column is equal to one.  Suppose $O = \left( o_{ij} \right) \in \R^{d
\times k}$ is nonnegative and column-stochastic.  Construct a discrete-time
process $\left\{ z_{t} \right\}$ to satisfy 
\[
 \mathbb{P} \left( z_{t}=i \mid x_{t} = j \right) = o_{ij}, 
\quad i \in \left\{ 1,\ldots,d \right\}, \;\;
j \in \left\{ 1,\ldots,k \right\}. 
\]
The matrix $O$ is called the observation matrix.  
If $\left\{ y_{t} \right\}$ and $\left\{ z_{t} \right\}$ have the same
law, we say $\left( O, Q \right)$ is an HMM realization of order $k$.  

When there exists a function 
$\phi \;:\; \left\{ 1,\ldots,k\right\} \rightarrow \left\{ 1,\ldots,d
\right\}$ such that $y_{t}=\phi(x_{t})$, then the observation matrix
can be selected as 
\begin{equation}
\label{eqn:observation_matrix}
 o_{ij} = 
 \left\{ 
 \begin{matrix}
  1 & \;\; \mbox{if} \;\;  \phi(j) = i, \\
  0 & \;\; \mbox{if} \;\;  \phi(j) \neq i.
 \end{matrix}
 \right.
\end{equation}
In this case, we say the observation is a deterministic function of the
state.  Note that each column of the matrix $O$ has exactly one nonzero
entry. 

The HMM-realization problem is to find a realization 
$\left( O,Q \right)$ given the probabilities 
$\mathbb{P}( y_{s}^{t} = u_{s}^{t})$ 
for any 
$s,t \in \Z$ and 
$u_{s}^{t} \in \left\{ 1,\ldots,d \right\}^{\left| t-s \right|+1}$.

\section{Reachable and Null Spaces}
\label{sec:reachable_and_null_spaces}

In this section, we review the reachable subspace and the null space for
an HMM introduced in \cite{ItoAmariKobayashi1992} 
and derive representations of these subspaces for the later discussion. 

Let $Q \in \R^{k \times k}$ be a matrix having a maximum modulus 
eigenvalue at one with the left eigenvector 
$e^{\trp} \in \R^{1 \times k}$ 
and the right eigenvector $\rho \in \R^{k}$ where $e$ is the vector
whose elements are all one.   
Let $\phi \;:\; 
\left\{ 1,\ldots,k\right\} \rightarrow \left\{ 1,\ldots,d
\right\}$ be a map and define 
$O = (o_{ij}) \in \R^{d \times k}$ by (\ref{eqn:observation_matrix}). 
Note that $Q$ needs not to be a nonnegative matrix but the sum of each
column is one.  
Let $I_{u} \in \R^{k \times k}$ ($u=1,\ldots,d$) 
be the diagonal matrix whose $(i,i)$th
element is one if $\phi(i)=u$ and zero otherwise.

We define the reachable subspace $V_{R}$ by
\begin{multline}
\label{eqn:reachable_subspace}
V_{R} = \sp 
\left\{ I_{u_{n}}Q
            \cdots I_{u_{1}}QI_{u_{0}} \rho \;:\; \right. \\ \left.
       u_{0}^{n}
       \in \left\{ 1,\ldots,d \right\}^{n+1}, 
       n=0,1,2,\ldots 
\right\}, 
\end{multline}
and the null space $V_{N}$ by 
\begin{equation}
\label{eqn:null_subspace}
V_{N} = 
\bigcap_{\substack{u_{0}^{n} 
      \in \left\{ 1,\ldots,d \right\}^{n+1},\\ 
      n=0,1,2,\ldots}}
\left\{ v \;:\; e^{\trp}I_{u_{n}}Q \cdots 
        I_{u_{1}}QI_{u_{0}}v = 0
\right\}.
\end{equation}
Define 
$V_{R,u}=I_{u}V_{R}$ ($u=1,\ldots,d$).  Then it follows that 
\begin{equation}
\label{eqn:reachbale_space_structure}
V_{R,u} \subset V_{R}, \quad 
V_{R} = \bigoplus_{u=1}^{d} V_{R,u}. 
\end{equation}
From the definitions (\ref{eqn:reachable_subspace}) and 
(\ref{eqn:null_subspace}), the following result is immediate. 
\begin{prop}
\label{prop:Q_I_invariant}
The reachable subspace $V_{R}$ is the smallest subspace which is
$Q$ as well as $I_{u}$-invariant ($u=1,\ldots,d$) and
 contains $\rho$. 
The null subspace $V_{N}$ is the largest subspace which is
$Q$ as well as $I_{u}$-invariant ($u=1,\ldots,d$) and
 contained in $\ker e^{\trp}$. 
\end{prop}

\begin{pf}
Since $\sum_{u=1}^{d} I_{u} = I$, $\rho = \sum_{u=1}^{d} I_{u}\rho \in
 V_{R}$. 
If $v \in V_{R}$, then 
$
 Qv = \sum_{u=1}^{d} I_{u}Qv \in V_{R}. 
$
Since $I_{u}I_{u}=I_{u}$ and $I_{u}I_{u^{\prime}}=0$ if
$u \neq u^{\prime}$, each $V_{R,u^{\prime}}$ is $I_{u}$
invariant and so is $V_{R}$. 
Conversely, if $V^{\prime}$ is $Q$-invariant as well as
 $I_{u}$-invariant and contains $\rho$, then 
$
I_{u_{n}}QI_{u_{n-1}}Q 
            \cdots I_{u_{1}}QI_{u_{0}} \rho \in V^{\prime}$ 
for any $(u_{0},\ldots,u_{n}) \in \left\{ 1,\ldots,d
 \right\}^{n+1}$ ($n=0,1,2,\ldots$). Since $V^{\prime}$ is a subspace,
 it contains any linear combination of this form.  Hence $V_{R} \subset
 V^{\prime}$. 
The proof for the null space is similar and is omitted. 
\qed
\end{pf}

\begin{assum}
\label{assumption:orthogonal_complement_1} 
The subspaces $V_{R,u}$ ($u=1,\ldots,d$) 
satisfy 
$V_{R,u} \not\subset \ker e^{\trp}$.
\end{assum}

Now, we give a basis of the reachable subspace.  

\begin{prop}
\label{prop:basis_reachable_subspace} 
Suppose Assumption \ref{assumption:orthogonal_complement_1} holds. 
Let $\hat{k}_{R} = \dim V_{R}$. 
There exist 
a full-column rank matrix 
$T_{R} = \left( t_{R,ij} \right) \in \R^{k \times \hat{k}_{R}}$ and 
a map  
$\hat{\phi}_{R} \;:\; 
\left\{ 1,\ldots,\hat{k}_{R} \right\} 
\rightarrow 
\left\{ 1,\ldots,d \right\}$ 
satisfying 
\begin{align}
\label{eqn:basis_reachable_subspace_1}
t_{R,ij} &= 0 \;\; \mbox{if} \;\; \phi(i) \neq \hat{\phi}_{R}(j), \\ 
\label{eqn:basis_reachable_subspace_2}
V_{R} &= \ran T_{R}, \\ 
\label{eqn:basis_reachable_subspace_3}
e^{\trp} T_{R} &= \hat{e}_{R}^{\trp}, 
\end{align}
where $\hat{e}_{R} \in \R^{\hat{k}_{R}}$ is the vector whose
elements are all one. 
\end{prop}

\begin{pf}
Choose a basis of $V_{R,u}$.  Then a basis of $V_{R}$ is constructed by
 collecting these bases.  
Define $\hat{\phi}_{R}$ by $\hat{\phi}_{R}(i)=j$ if 
$\sum_{u=1}^{j-1} \hat{k}_{u} < i \leq \sum_{u=1}^{j} \hat{k}_{u}$. 
Let $\hat{k}_{R}=\sum_{u=1}^{d} \hat{k}_{u} = \dim V_{R}$.  
Arrange the basis to form a matrix 
$T_{R} \in \R^{k \times \hat{k}_{R}}$ so that 
 the $i$th columns of $T_{R}$ ($i \in \hat{\phi}_{R}^{-1}(u)$) are a basis of
 $V_{R,u}$. 
From the construction, $T_{R}$ 
satisfies (\ref{eqn:basis_reachable_subspace_1}). 
From Assumption \ref{assumption:orthogonal_complement_1},
 we can select those columns to satisfy 
$e^{\trp}T=\hat{e}^{\trp}$ by 
perturbing a little bit and rescaling.  
\qed
\end{pf}

\begin{prop}
\label{prop:basis_dual_null_subspace}
Let 
$\hat{k}_{N} = k - \dim V_{N}$. 
There exist 
a full-row rank matrix 
$T_{N} = \left( t_{N,ij} \right) \in \R^{\hat{k}_{N} \times k}$ and 
a map 
$\hat{\phi}_{N} \;:\; 
\left\{ 1,\ldots,\hat{k}_{N} \right\} 
\rightarrow 
\left\{ 1,\ldots,d \right\}$ 
satisfying 
\begin{align}
\label{eqn:basis_dual_null_subspace_1}
t_{N,ij} &= 0 
\;\; \mbox{if} \;\; \hat{\phi}_{N}(i) \neq \phi(j), \\
\label{eqn:basis_dual_null_subspace_2}
V_{N} &= \ker T_{N}, \\ 
\label{eqn:basis_dual_null_subspace_3}
e^{\trp} &= \hat{e}_{N}^{\trp}T_{N}, 
\end{align}
where $\hat{e}_{N} \in \R^{\hat{k}_{N}}$ is the vector whose
elements are all one. 
\end{prop}

\begin{rem}
\label{remark:basis_reachable_space}
Note that the $Q$-invariance of $V_{R}$ 
and (\ref{eqn:basis_reachable_subspace_2}) imply that there exists 
$\hat{Q}_{R} \in \R^{\hat{k}_{R} \times \hat{k}_{R}}$ such that 
$QT_{R}=T_{R}\hat{Q}_{R}$. 
Since $\rho \in V_{R}$, there exists $\hat{\rho}_{R} \in \R^{\hat{k}_{R}}$
 such that $\rho = T_{R}\hat{\rho}_{R}$.  
From (\ref{eqn:basis_reachable_subspace_3}) and the definition of
 $\hat{\rho}_{R}$, we have 
\begin{align*}
\hat{e}_{R}^{\trp}\hat{Q}_{R}
&=e^{\trp}T_{R}\hat{Q}_{R}
=e^{\trp}QT_{R} = e^{\trp}T_{R} = \hat{e}_{R}^{\trp}, \\
T_{R}\hat{Q}_{R} \hat{\rho}_{R} &= 
QT_{R} \hat{\rho}_{R} = Q \rho = \rho \;\; 
\Rightarrow  \;\;
\hat{Q}_{R}\hat{\rho}_{R}=\hat{\rho}_{R}, 
\end{align*}
which means 
$\hat{Q}_{R}$ has an eigenvalue at one with the left eigenvector
 $\hat{e}_{R}$ and the right eigenvector $\hat{\rho}_{R}$. 
Similarly, 
the $Q$-invariance of $V_{N}$ 
and (\ref{eqn:basis_dual_null_subspace_2}) 
imply that there exists 
$\hat{Q}_{N} \in \R^{\hat{k}_{N} \times \hat{k}_{N}}$ such that 
$\hat{Q}_{N}T_{N}=T_{N}Q$. Define $\hat{\rho}_{N}=T_{N}\rho$.  
Then $\hat{Q}_{N}$ has an eigenvalue at one with the left eigenvector
 $\hat{e}_{N}$ and the right eigenvector $\hat{\rho}_{N}$. 
\end{rem}

\section{Tensor Decomposition}
\label{sec:tensor_decomposition}

In \cite{HuangGeKakadeDahleh2016}, 
a third-order tensor was introduced to discuss the HMM realization
problem. Tensor decomposition methods are exploited to solve the minimal
HMM realization problem and showed that the minimal order of HMM
realization is equal to the rank of the tensor excluding a measure zero
set of parameter space.  

In this section, we concentrate on the case where the observation is
deterministic and may not be captured by the generic argument.  
Define a third-order tensor 
$M \in \R^{d^{n} \times d^{n} \times d}$ 
whose 
$\left( L(u_{1}^{n}),L(u_{-1}^{-n}),u_{0} \right)$th element is 
given by 
\begin{equation}
\label{eqn:tensor_M}
 M_{L(u_{1}^{n}),L(u_{-1}^{-n}),u_{0}} 
= \mathbb{P} \left\{ y_{-n}^{n} = u_{-n}^{n} \right\}, 
\end{equation}
where $L$ is the index map $\left\{ 1,\ldots,d \right\}^{n} \rightarrow
\Z$ that assigns an element of
$\left\{ 1,\ldots,d \right\}^{n}$ 
to its $d$-digit number, or namely 
\[
 L(u_{1}^{n}) 
= \left( u_{1}-1 \right)d^{n-1}+\left( u_{2}-1 \right)d^{n-2}+
  \cdots+u_{n} 
\]
for $u_{1}^{n} = \left( u_{1}, \ldots, u_{n} \right)$.  

Suppose $\left( O, Q \right)$ is an HMM realization of order $k$
generating the process $\left\{ y_{t} \right\}$. 
Then the tensor $M$ can be decomposed into a sum of
rank-one tensors 
\begin{equation}
\label{eqn:decomposition_M} 
 M = A \otimes B \otimes C, 
\end{equation}
where $A \in \R^{d^{n} \times k}$, 
$B \in \R^{d^{n} \times k}$, and 
$C \in \R^{d \times k}$ are given by 
\begin{align}
\label{eqn:A_row_component_new}
A_{L(u_{1}^{n})} &=
e^{\trp}
I_{u_{n}}Q
\cdots
I_{u_{1}}Q, \\
\label{eqn:B_row_component_new}
B_{L(u_{1}^{n})} &=
\rho^{\trp}
I_{u_{n}}Q^{\trp}
\cdots
I_{u_{1}}Q^{\trp}, \\
\label{eqn:C_new}
C &= O, 
\end{align}
where 
$A_{L(u_{1}^{n})}$ and 
$B_{L(u_{1}^{n})}$ denote the 
$L(u_{1}^{n})$th row of the matrices 
$A$ and $B$, respectively.  
Note that 
(\ref{eqn:A_row_component_new}) and (\ref{eqn:B_row_component_new}) 
are the row-wise expression of the recursive formula using Khatri-Rao
products.  Note also that the definition of $B$ and $C$ in 
(\ref{eqn:B_row_component_new}) and (\ref{eqn:C_new}) is modified from 
\cite{HuangGeKakadeDahleh2016} by scalar multiplication of the columns.

\subsection{Reduction using reachable subspace}

If the reachable space $V_{R}$ is not the whole space $\R^{k}$, then 
the number of the rank-one tensors in (\ref{eqn:decomposition_M}) can be
reduced to $\hat{k}_{R} = \dim V_{R}$. 

\begin{thm}
\label{thm:reduction_reachable_subspace}
Suppose $\hat{k}_{R}=\dim V_{R}<k$.  
Define $T_{R} \in \R^{k \times \hat{k}_{R}}$ and 
$\hat{\phi}_{R} \;:\; 
\left\{ 1,\ldots,\hat{k}_{R} \right\} 
\rightarrow 
\left\{ 1,\ldots,d \right\}$ 
as in Proposition \ref{prop:basis_reachable_subspace}. 
Define $\hat{Q}_{R} \in \R^{\hat{k}_{R} \times \hat{k}_{R}}$, 
$\hat{O}_{R} \in \R^{d \times \hat{k}_{R}}$, 
$\hat{e}_{R}$ , and 
$\hat{\rho}_{R}$ by 
$QT_{R}=T_{R}\hat{Q}_{R}$, 
$\hat{O}_{R}=OT_{R}$, 
$\hat{e}_{R}^{\trp}=e^{\trp}T_{R}$, and
$\rho = T_{R}\hat{\rho}_{R}$, 
respectively. 
Define 
$\hat{A}_{R} \in \R^{d^{n} \times \hat{k}_{R}}$, 
$\hat{B}_{R} \in \R^{d^{n} \times \hat{k}_{R}}$, and  
$\hat{C}_{R} \in \R^{d \times \hat{k}_{R}}$ by 
\begin{align}
\label{eqn:A_row_component_reduced_r}
\hat{A}_{R,L(u_{1}^{n})} &=
\hat{e}_{R}^{\trp}
\hat{I}_{R,u_{n}}\hat{Q}_{R}
\cdots
\hat{I}_{u_{1}}\hat{Q}_{R}, \\
\label{eqn:B_row_component_reduced_r}
\hat{B}_{R,L(u_{1}^{n})} &=
\hat{\rho}_{R}^{\trp}
\hat{I}_{R,u_{n}}\hat{Q}_{R}^{\trp}
\cdots
\hat{I}_{R,u_{1}}\hat{Q}_{R}^{\trp}, \\
\label{eqn:C_reduced_r}
\hat{C}_{R} &= \hat{O}_{R}, 
\end{align}
where $\hat{I}_{R,u}$ is the diagonal matrix whose $\left( i,i
 \right)$th element is one if $\hat{\phi}_{R}(i)=u$ and zero otherwise. 
Then 
\[
 M = A \otimes B \otimes C 
= \hat{A}_{R} \otimes \hat{B}_{R} \otimes \hat{C}_{R}
\]
holds.
\end{thm}

\begin{pf}
From $QT_{R}=T_{R}\hat{Q}_{R}$ and 
$I_{u}T_{R}=T_{R}\hat{I}_{R,u}$, we have
\begin{align*}
A_{L(u_{1}^{n})}T_{R} 
&=
e^{\trp}
I_{u_{n}}Q
\cdots
I_{u_{1}}Q T_{R} \\
&=
e^{\trp}
I_{u_{n}}Q
\cdots
I_{u_{1}} T_{R} \hat{Q}_{R}\\
&=
e^{\trp}
I_{u_{n}}Q
\cdots
T_{R} \hat{I}_{R,u_{1}} \hat{Q}_{R}\\
& \vdots \\
&=
e^{\trp}T_{R}
\hat{I}_{R,u_{n}}\hat{Q}_{R}
\cdots
\hat{I}_{R,u_{1}} \hat{Q}_{R} \\
&=
\hat{e}_{R}^{\trp}
\hat{I}_{R,u_{n}}\hat{Q}_{R}
\cdots
\hat{I}_{R,u_{1}} \hat{Q}_{R} 
 = 
\hat{A}_{R,L(u_{1}^{n})} \\
\hat{B}_{R,L(u_{1}^{n})}T_{R}^{\trp} 
&=
\hat{\rho}_{R}^{\trp}
\hat{I}_{R,u_{n}}\hat{Q}_{R}^{\trp}
\cdots
\hat{I}_{R,u_{1}}\hat{Q}_{R}^{\trp}T_{R}^{\trp} \\
&=
\hat{\rho}_{R}^{\trp}
\hat{I}_{R,u_{n}}\hat{Q}_{R}^{\trp}
\cdots
\hat{I}_{R,u_{1}}T_{R}^{\trp}Q^{\trp} \\
&=
\hat{\rho}_{R}^{\trp}
\hat{I}_{R,u_{n}}\hat{Q}_{R}^{\trp}
\cdots
T_{R}^{\trp}I_{u_{1}}Q^{\trp} \\
&\vdots \\
&=
\hat{\rho}_{R}^{\trp}T_{R}^{\trp}
I_{u_{n}}Q^{\trp}
\cdots
I_{u_{1}}Q^{\trp} \\
&=
\rho^{\trp}
I_{u_{n}}Q^{\trp}
\cdots
I_{u_{1}}Q^{\trp} 
=
B_{L(u_{1}^{n})}. 
\end{align*}
Hence
\begin{equation}
\label{eqn:relations_c} 
 AT_{R}=\hat{A}_{R}, \quad
 B = \hat{B}_{R}T_{R}^{\trp}, \quad
 CT_{R} = \hat{C}_{R}.
\end{equation}
Denote the column vectors of 
$\hat{A}_{R}$, 
$\hat{B}_{R}$, and 
$\hat{C}_{R}$ by 
$\hat{a}_{i}$, 
$\hat{b}_{i}$, and 
$\hat{c}_{i}$ ($i=1,\ldots,\hat{k}$), respectively, and 
the column vectors of 
$A$, $B$, and $C$ by 
$a_{i}$, $b_{i}$, and $c_{i}$ ($i=1,\ldots,k$), respectively. 
From (\ref{eqn:relations_c}), we have
\begin{align*}
A \otimes B \otimes C
&=
\sum_{i=1}^{k} a_{i} \otimes b_{i} \otimes c_{i} \\
&= 
\sum_{i=1}^{k} 
  a_{i} \otimes 
  \left( \sum_{j=1}^{\hat{k}_{R}} t_{R,ij}\hat{b}_{j} \right) 
  \otimes c_{i} \\
&= 
\sum_{j=1}^{\hat{k}_{R}}
\sum_{i=1}^{k} 
  t_{R,ij} a_{i} \otimes \hat{b}_{j} \otimes c_{i} \\
&=
\sum_{j=1}^{\hat{k}_{R}}
\sum_{i \in \phi^{-1}(\hat{\phi}_{R}(j))} 
  t_{R,ij} a_{i} \otimes \hat{b}_{j} \otimes c_{i} \\
&=
\sum_{j=1}^{\hat{k}_{R}}
\sum_{i \in \phi^{-1}(\hat{\phi}_{R}(j))} 
  t_{R,ij} a_{i} \otimes \hat{b}_{j} \otimes \hat{c}_{j} \\
&=
\sum_{j=1}^{\hat{k}_{R}}
\left( 
\sum_{i=1}^{k}
  t_{R,ij} a_{i} \right) 
\otimes \hat{b}_{j} \otimes \hat{c}_{j}  \\
&=
\sum_{j=1}^{\hat{k}_{R}}
  \hat{a}_{j} \otimes \hat{b}_{j} \otimes \hat{c}_{j} = 
\hat{A}_{R} \otimes \hat{B}_{R} \otimes \hat{C}_{R}.
\end{align*}
\qed
\end{pf}

\subsection{Reduction using null space}

If the null space $V_{N}$ is not the zero subspace, then 
the number of the rank-one tensors in (\ref{eqn:decomposition_M}) can be
reduced to $\hat{k}_{N} = k-\dim V_{N}$.

\begin{thm}
\label{thm:reduction_null_subspace}
Suppose $\hat{k}_{N} = k-\dim V_{N}<k$.  
Define $T_{N} \in \R^{\hat{k}_{N} \times k}$ and 
$\hat{\phi}_{N} \;:\; 
\left\{ 1,\ldots,\hat{k}_{N} \right\} 
\rightarrow 
\left\{ 1,\ldots,d \right\}$ 
as in Proposition \ref{prop:basis_dual_null_subspace}. 
Define $\hat{Q}_{N} \in \R^{\hat{k}_{N} \times \hat{k}_{N}}$, 
$\hat{O}_{N} \in \R^{d \times \hat{k}_{N}}$, 
$\hat{e}_{N}$, and 
$\hat{\rho}_{N}$ by
$\hat{Q}_{N}T_{N}=T_{N}Q$, 
$\hat{O}_{N}T_{N}=O$, 
$\hat{e}_{N}^{\trp}T_{N}=e^{\trp}$, and
$\hat{\rho}_{N}=T_{N}\rho$, 
respectively. 
Define 
$\hat{A}_{N} \in \R^{d^{n} \times \hat{k}_{N}}$, 
$\hat{B}_{N} \in \R^{d^{n} \times \hat{k}_{N}}$, and  
$\hat{C}_{N} \in \R^{d \times \hat{k}_{N}}$ by 
\begin{align}
\label{eqn:A_row_component_reduced_n}
\hat{A}_{N,L(u_{1}^{n})} &=
\hat{e}_{N}^{\trp}
\hat{I}_{N,u_{n}}\hat{Q}_{N}
\cdots
\hat{I}_{N,u_{1}}\hat{Q}_{N}, \\
\label{eqn:B_row_component_reduced_n}
\hat{B}_{N,L(u_{1}^{n})} &=
\hat{\rho}_{N}^{\trp}
\hat{I}_{N,u_{n}}\hat{Q}_{N}^{\trp}
\cdots
\hat{I}_{N,u_{1}}\hat{Q}_{N}^{\trp}, \\
\label{eqn:C_reduced_n}
\hat{C}_{N} &= \hat{O}_{N}, 
\end{align}
where $\hat{I}_{N,u}$ is the diagonal matrix whose $\left( i,i
 \right)$th element is one if $\hat{\phi}_{N}(i)=u$ and zero otherwise. 
Then 
\[
 M = A \otimes B \otimes C 
= \hat{A}_{N} \otimes \hat{B}_{N} \otimes \hat{C}_{N}
\]
holds.
\end{thm}

\begin{pf}
From $\hat{Q}_{N}T_{N}=T_{N}Q$ and 
$\hat{I}_{N,u}T_{N}=T_{N}I_{u}$, we have
\[
 \hat{A}_{N}T_{N}=A, \quad
 \hat{B}_{N} = BT_{N}^{\trp}, \quad
 \hat{C}_{N}T_{N}=C.
\]
The rest of the proof is similar to the proof of 
Theorem \ref{thm:reduction_reachable_subspace} and is omitted. 
\qed
\end{pf}

\subsection{Reduction using effective space}

When the reachable subspace is not the whole space and the null space is
not the zero subspace, we may have the situation where both Theorems 
\ref{thm:reduction_reachable_subspace} and 
\ref{thm:reduction_null_subspace} can be applied.  
In fact, by using the notion of the effective space in 
\cite{ItoAmariKobayashi1992}, we can reduce the number of rank-one
tensors in the decomposition (\ref{eqn:decomposition_M}). 

The effective space is defined as 
\[
\left(  V_{R}+V_{N} \right)/V_{N},  
\]
and its dimension is given by
\begin{align*}
\hat{k}&= 
\dim \left(  V_{R}+V_{N} \right)/V_{N} \\
&= \dim \left(  V_{R}+V_{N} \right) - \dim V_{N} \\
&= \dim V_{R} - \dim \left( V_{R} \cap V_{N} \right). 
\end{align*}

\begin{prop}
\label{prop:basis_effective_space}
Select $T_{R}$ and $T_{N}$ as in Propositions 
\ref{prop:basis_reachable_subspace} and 
\ref{prop:basis_dual_null_subspace}, respectively.  
Then $\hat{k}=\rank T_{N}T_{R}$.  
There exist a full row rank matrix  
$T = (t_{ij}) \in \R^{\hat{k} \times \hat{k}_{r}}$ 
and a map 
$\hat{\phi} \;:\; 
\left\{ 1,\ldots,\hat{k} \right\} \rightarrow
\left\{ 1,\ldots,d \right\}$
satisfying 
\begin{align}
\label{eqn:basis_effective_space_1}
 t_{ij} &= 0 \; \mbox{if}\; \hat{\phi}(i) \neq \hat{\phi}_{R}(j), \\
\label{eqn:basis_effective_space_2}
 \ker T &= \ker T_{N}T_{R}, \\
\label{eqn:basis_effective_space_3}
 \hat{e}_{R}^{\trp} &= \hat{e}^{\trp}T, \\ 
\label{eqn:basis_effective_space_4}
 \hat{Q}_{R} \ker T & \subset \ker T.  
\end{align}
\end{prop}

\begin{pf}
Because $T_{N}$ and $T_{R}$ have the block structure, 
we can select $\hat{k}$ linearly independent rows of $T_{N}T_{R}$ to
construct $T = (t_{ij}) \in \R^{\hat{k} \times \hat{k}_{r}}$ satisfying 
(\ref{eqn:basis_effective_space_1}) and 
(\ref{eqn:basis_effective_space_2}). 
The condition (\ref{eqn:basis_effective_space_3}) 
is satisfied by scaling and perturbation if necessary.  
From $T_{N}T_{R}\hat{Q}_{R}=T_{N}QT_{R}=\hat{Q}_{N}T_{N}T_{R}$, 
$\ker T = \ker T_{N}T_{R}$ is $\hat{Q}_{R}$ invariant.   
\qed
\end{pf}

\begin{thm}
\label{thm:reduction_effective_subspace}
Define 
$T \in \R^{\hat{k} \times \hat{k}_{r}}$ and 
$\hat{\phi} \;:\; 
\left\{ 1,\ldots,\hat{k} \right\} \rightarrow
\left\{ 1,\ldots,d \right\}$ as in 
Proposition \ref{prop:basis_effective_space}. 
Define $\hat{Q} \in \R^{\hat{k} \times \hat{k}}$, 
$\hat{O} \in \R^{d \times \hat{k}}$, 
$\hat{e}$, and 
$\hat{\rho}$ by
$\hat{Q}T=T\hat{Q}_{R}$, 
$\hat{O}T=\hat{O}_{R}$, 
$\hat{e}^{\trp}T=\hat{e}_{R}^{\trp}$, and
$\hat{\rho}=T\hat{\rho}_{R}$, 
respectively. 
Define 
$\hat{A} \in \R^{d^{n} \times \hat{k}}$, 
$\hat{B} \in \R^{d^{n} \times \hat{k}}$, and  
$\hat{C} \in \R^{d \times \hat{k}}$ by 
\begin{align}
\label{eqn:A_row_component_reduced_e}
\hat{A}_{L(u_{1}^{n})}
&=
\hat{e}^{\trp}
\hat{I}_{u_{n}}\hat{Q}
\cdots
\hat{I}_{u_{1}}\hat{Q}, \\
\label{eqn:B_row_component_reduced_e}
\hat{B}_{L(u_{1}^{n})} &=
\hat{\rho}^{\trp}
\hat{I}_{u_{n}}\hat{Q}^{\trp}
\cdots
\hat{I}_{u_{1}}\hat{Q}^{\trp}, \\
\label{eqn:C_reduced_e}
\hat{C} &= \hat{O}, 
\end{align}
where $\hat{I}_{u}$ is the diagonal matrix whose $\left( i,i
 \right)$th element is one if $\hat{\phi}(i)=u$ and zero otherwise. 
Then 
\[
 M = A \otimes B \otimes C 
= \hat{A} \otimes \hat{B} \otimes \hat{C} 
\]
holds.
\end{thm}

\begin{pf}
We apply Theorem \ref{thm:reduction_null_subspace} to the system $\left(
 \hat{O}_{R}, \hat{Q}_{R} \right)$.  
\qed
\end{pf}

\begin{rem}
\label{remark:generalized_Hankel_matrix}
The generalized Hankel matrix introduced in \cite{Picci1978} 
takes the following form (see also \cite{Anderson1999}). 
\[
H = 
\begin{bmatrix}
 H^{(00)} &   H^{(01)} &   H^{(02)} & \cdots \\
 H^{(10)} &   H^{(11)} &   H^{(12)} & \vdots \\
 H^{(20)} &   H^{(21)} &   H^{(22)} & \vdots \\
 \vdots & \vdots & \vdots & \ddots
\end{bmatrix},
\]
where $H^{(ij)}$ is a $d^{i} \times d^{j}$ matrix whose 
$\left( L(u_{0}^{-(i-1)}), L(u_{1}^{j}) \right)$th element is given by 
\[
 H^{(ij)}_{L(u_{0}^{-(i-1)}), L(u_{1}^{j})} = 
\mathbb{P} \left\{ y_{-(i-1)}^{j} = u_{-(i-1)}^{j} \right\}. 
\]
Let $\left( O, Q \right)$ be a realization of the HMM of order $k$. 
Define 
\begin{align*}
 \Theta &= \begin{bmatrix}
	    \Theta^{(0)} & \Theta^{(1)} & \Theta^{(2)} & \cdots 
	   \end{bmatrix}, \\
 \Gamma &= \begin{bmatrix}
	    \Gamma^{(0)} & \Gamma^{(1)} & \Gamma^{(2)} & \cdots 
	   \end{bmatrix}, 
\end{align*}
where 
$\Theta^{(i)}$ is a $k \times d^{i}$ matrix whose 
$L(u_{0}^{-(i-1)})$th column is given by 
\[
\Theta^{(i)}_{L(u_{0}^{-(i-1)})}
= I_{u_{0}}Q \cdots I_{u_{-(i-1)}}\rho, 
\]
and 
$\Gamma^{(i)}$ is a $k \times d^{i}$ matrix whose 
$L(u_{1}^{i})$th column is given by 
\[
 \Gamma^{(i)}_{L(u_{1}^{i})}
= I_{u_{i}}Q^{\trp} \cdots I_{u_{1}} e. 
\]
Then, we have
\[
 H = \Theta^{\trp}\Gamma
\]
holds, and thus $\rank H$ is at most $k$.  
Notice that the columns of $\Theta$ generate the reachable subspace 
(\ref{eqn:reachable_subspace}) and 
the columns of $\Gamma$ generate the orthogonal complement of 
the null space (\ref{eqn:null_subspace}).  
Hence, $\rank H$ is equal to the dimension of the effective subspace.
\end{rem}

\subsection{Example}

This example is modified from the example discussed in
\cite{Vidyasagar2011} (originally in 
\cite{FoxRubin1968} and 
\cite{DharmadhikariNadkarni1970}). 
Let $\lambda \in \left( 0, 0.5 \right]$ and 
$\alpha = 2\pi/m$ for some $m \in \left\{ 3,4,\ldots \right\}$. 
Let $\zeta = e^{j\alpha}$. 
Note that the example was intended to illustrate the case where the
generalized Hankel matrix has finite rank but it does not have a
finite-state Markov model realization; in this case, $\alpha$ is
selected to be non-commensurate to $\pi$.  The example in this section
assumes that $\alpha$ and $\pi$ are commensurate.

Suppose that a stationary discrete-time random process 
$\left\{ y_{t} \right\}$ taking binary values $\left\{ 1,2 \right\}$ 
has an HMM realization whose 
transition matrix $Q \in \R^{(m+1) \times (m+1)}$ 
and observation matrix $O \in \R^{2 \times (m+1)}$ are given by 
\begin{align*}
  Q &= 
 \begin{bmatrix}
  1-\sum_{i=1}^{m-1} 
   \frac{\lambda^{i}}{1-\lambda^{m}} \sin^{2} \frac{i\pi}{m} & 
  1-\lambda^{m} & 0 & \cdots & 0 \\
  \frac{\lambda}{1-\lambda^{m}} \sin^{2} \frac{\pi}{m} & 0 & 1 
  & \vdots & 0 \\
  \vdots & \vdots & \vdots & \ddots & \vdots \\
   \frac{\lambda^{m-1}}{1-\lambda^{m}} \sin^{2} \frac{(m-1)\pi}{m} & 
   0 & 0 & \vdots & 1 \\
   0 & \lambda^{m} & 0 & \cdots & 0
 \end{bmatrix}, \\
 O &= 
 \begin{bmatrix}
  1 & 0 & 0 & \cdots & 0 \\
  0 & 1 & 1 & \cdots & 1
 \end{bmatrix}. 
\end{align*}
Let 
\begin{align*}
p_{0} 
&= 
\begin{bmatrix}
  1 \\
  0 \\
  0 \\
 \vdots \\
  0
\end{bmatrix}, \quad
p_{1}
=
\begin{bmatrix}
 0 \\
 1 \\
 \lambda \\
 \vdots \\
 \lambda^{m-1}
\end{bmatrix}, \\
p_{2}
&=
\begin{bmatrix}
 0 \\
 1 \\
 \lambda \zeta \\
 \vdots \\
 \lambda^{m-1} \zeta^{m-1}
\end{bmatrix}, \quad
p_{3}
=
\begin{bmatrix}
 0 \\
 1 \\
 \lambda \zeta^{-1} \\
 \vdots \\
 \lambda^{m-1} \zeta^{-(m-1)}
\end{bmatrix}.
\end{align*}
Then we can show that 
$\left\{ p_{0}, p_{1}, \left( p_{2}+p_{3} \right)/2, 
 \left( p_{2}-p_{3} \right)/2j \right\}$ is a basis of the reachable
 space $V_{R}$.  By scaling, 
\begin{align*}
 T_{R} &= 
 \begin{bmatrix}
   p_{0} & p_{1} & p_{2} & p_{3}
 \end{bmatrix}
 \begin{bmatrix}
   1 & 0 & 0 & 0 \\
   0 & 1 & 0 & 0 \\
   0 & 0 & \frac{1}{2} & \frac{1}{2j} \\
   0 & 0 & \frac{1}{2} & -\frac{1}{2j}
 \end{bmatrix}
\\ &  \times 
 \begin{bmatrix}
   1 & 0 & 0 & 0 \\
   0 & \frac{1-\lambda^{m}}{1-\lambda} & 0 & 0 \\
   0 & 0 & \frac{(1-\lambda^{m})(1-\lambda \cos \alpha)}
 {1-2\lambda \cos \alpha + \lambda^{2}} & 0 \\
   0 & 0 & 0 & 
 \frac{\lambda (1-\lambda^{m}) \sin \alpha}
 {1-2\lambda \cos \alpha + \lambda^{2}}
 \end{bmatrix}^{-1}
\end{align*}
satisfies 
(\ref{eqn:basis_reachable_subspace_1}), 
(\ref{eqn:basis_reachable_subspace_2}), and 
(\ref{eqn:basis_reachable_subspace_3}) 
in Proposition \ref{prop:basis_reachable_subspace}. 
Then the matrices $\hat{Q}_{R} \in \R^{4 \times 4}$ and 
$\hat{O} \in \R^{2 \times 4}$ defined by 
$QT_{R}=T_{R}\hat{Q}_{R}$ and 
$\hat{O}_{R}=OT_{R}$ are given by 
\begin{align*}
 \hat{Q}_{R} 
&=
\left[ 
 \begin{matrix}
  \eta & 1-\lambda \\ 
  \frac{\lambda}{2(1-\lambda)} & \lambda \\ 
  \frac{\lambda \cos \alpha (\lambda \cos \alpha -1)}
       {2(1-2\lambda \cos \alpha + \lambda^{2})} & 0 \\
  \frac{\lambda^{2} \sin^{2} \alpha}
       {2(1-2\lambda \cos \alpha + \lambda^{2})} & 0 
 \end{matrix} \right. 
\\ & \hspace{6em} \left. 
 \begin{matrix}
  \frac{1-2\lambda \cos \alpha +\lambda^{2}}{1-\lambda \cos \alpha} & 0 \\
  0 & 0 \\
  \lambda \cos \alpha & 1-\lambda \cos \alpha \\
  -\frac{\lambda^{2} \sin^{2} \alpha}{1-\lambda \cos \alpha} & 
  \lambda \cos \alpha
 \end{matrix}
\right], \\
\hat{O}_{R} &=
 OT_{R} = 
 \begin{bmatrix}
   1 & 0 & 0 & 0 \\
   0 & 1 & 1 & 1
 \end{bmatrix}, 
\end{align*}
where 
\[
\eta = 
2-3\lambda -3\lambda\cos\alpha 
+ \lambda^{2}+5\lambda^{2}\cos\alpha-2\lambda^{3}.
\]
Let 
\begin{align*}
  \hat{I}_{R,1} &=
 \begin{bmatrix}
   1 & 0 & 0 & 0 \\
   0 & 0 & 0 & 0 \\
   0 & 0 & 0 & 0 \\
   0 & 0 & 0 & 0 
 \end{bmatrix}, \quad
 \hat{I}_{R,2} = 
 \begin{bmatrix}
   0 & 0 & 0 & 0 \\
   0 & 1 & 0 & 0 \\
   0 & 0 & 1 & 0 \\
   0 & 0 & 0 & 1
 \end{bmatrix}, \\
\hat{e}_{R} &=
 \begin{bmatrix}
  1 \\
  1 \\
  1 \\
  1
 \end{bmatrix}, \quad
\hat{\rho}_{R} = 
\begin{bmatrix}
 \rho_{1} \\
 \frac{\rho_{1}\lambda}{2(1-\lambda)^{2}} \\
 \frac{\rho_{1}\lambda(\cos\alpha-\lambda)(\lambda\cos\alpha-1)}
      {2(1-2\lambda\cos\alpha + \lambda^{2})^{2}} \\
  \frac{\rho_{1}\lambda^{2}\sin^{2}\alpha}
       {2(1-2\lambda\cos\alpha + \lambda^{2})^{2}} 
\end{bmatrix},
\end{align*}
where $\rho_{1}$ is selected to satisfy
$\hat{e}_{R}^{\trp}\hat{\rho}_{R}=1$. 
Define $\hat{A}_{R} \in \R^{2^{n} \times 4}$, 
$\hat{B}_{R} \in \R^{2^{n} \times 4}$, and 
$\hat{C}_{R} \in \R^{2 \times 4}$ by 
(\ref{eqn:A_row_component_reduced_r}), 
(\ref{eqn:B_row_component_reduced_r}), and 
(\ref{eqn:C_reduced_r}), respectively.  
Then $\hat{A}_{R} \otimes \hat{B}_{R} \otimes \hat{C}_{R}$ is 
a sum of four rank-one tensors for the 
third-order tensor (\ref{eqn:tensor_M}) which is originally 
written as a sum of $(m+1)$ rank-one tensors.  
So, if $m \geq 4$, then we have a reduced number of rank-one tensors for 
the decomposition (\ref{eqn:decomposition_M}).

\section{Conclusion}
\label{sec:conclusion}

This paper considered the minimum HMM realization problem using tensor
decomposition methods. If the observation is deterministic, or the
Kruskal rank of the observation matrix equals one, then the third-order
tensor can be decomposed to a sum of rank-one tensors whose number is
not greater than the dimension of the effective space.  Since the
dimension of the effective space is equal to the rank of the generalized
Hankel matrix, the tensor decomposition ends up giving a minimal
quasi-realization. This means that determining the minimum number of
states to realize a stationary stochastic process with a finite alphabet
for the class of HMM's with deterministic observation is yet
unresolved.


\end{document}